\documentclass[11pt]{amsart}
\usepackage{amsfonts,amsmath,amssymb}

\newtheorem{theorem}{Theorem}[section]
\newtheorem{lemma}[theorem]{Lemma}
\newtheorem{proposition}[theorem]{Proposition}

\newtheorem{definition}{Definition}[section]
\newtheorem{corollary}[theorem]{Corollary}
\newtheorem{remark}[theorem]{Remark}

\newcommand{\cl}[1]{\mathcal{#1}} 

\newcommand{\Map}[1]{\mathrm{Map}(#1)} 
\newcommand{\Ref}[1]{\mathrm{Ref}(#1)} 
\newcommand{\Lat}[1]{\mathrm{Lat}(#1)} 
\newcommand{\Alg}[1]{\mathrm{Alg}(#1)} 

\title{TRO EQUIVALENT ALGEBRAS}
\author{G.K. ELEFTHERAKIS}

\address{GEORGE K. ELEFTHERAKIS\\
 Department of Mathematics,\\ University of
Athens, \\ Panepistimioupolis 157 84\\ Athens Greece}
\email{gelefth@math.uoa.gr}

\thanks{This
research was partly supported by Special Account Research Grant
No. 70/3/7463 of the University of Athens.}



\begin{document}

\maketitle

\begin{abstract}
In this work we study a new equivalence relation between $w^*$ closed algebras of operators
on Hilbert spaces. The algebras $\cl{A}$ and $\cl{B}$ are called \emph{TRO equivalent} if
 there exists a ternary ring of operators $\cl{M}$ (i.e. $\cl{M}\cl{M}^*\cl{M}\subset \cl{M}$)
 such that $\cl{A}= span( \cl{M}^*\cl{B}\cl{M})^{-w^*}$ and $\cl{B}=span( \cl{M}\cl{A}\cl{M}^*)^{-w^*}.$
 We prove that two reflexive algebras are TRO equivalent if and only if  there exists a $*$
isomorphism between the commutants of their diagonals mapping the invariant projection
 lattice of the first algebra onto the lattice of the second one. 
\end{abstract}

\section{Introduction}

A linear space $\cl{M}$ of operators between two Hilbert spaces satisfying 
$$\cl{M}\cl{M}^*\cl{M}\subset \cl{M}$$ 
is called a \textbf{ ternary ring of operators  (TRO)}.

TRO's were introduced in  \cite{hes} and constitute a generalisation of
selfadjoint operator algebras \cite{har}, \cite{zet}.  They have many properties similar to
$C^*$-algebras and von Neumann algebras. Recently, these objects have been studied
from the point of view of operator space theory, in which they 
play an important role \cite{eor}, \cite{ruankaur}, \cite{ruan}.

In \cite{kt}, TRO's were studied from a different angle, namely as {\em normalisers}
of operator algebras:

If $\cl{A}\subset B(H)$ and $\cl{B}\subset B(K)$ are $w^*$ closed operator algebras, not necessarily
selfadjoint, an operator $T\in B(H,K)$ is said to {\em normalise the algebra $\cl{A}$ into $\cl{B}$} if
$T^*\cl{B}T\subset \cl{A}$ and  $T\cl{A}T^*\subset \cl{B}.$ It is shown in \cite{kt} that
such a normaliser $T$ defines a TRO $\cl{M}_T$ consisting of normalisers from $\cl{A}$ into $\cl{B}:$
$$ \cl{M}_T^* \cl{B}\cl{M}_T\subset \cl{A}  \;\;\text{and }\;\;\cl{M}_T \cl{A}\cl{M}_T^*\subset \cl{B}.$$

In the present paper we are interested in a stronger situation, namely in the
existence of a TRO $\cl{M}$ so that
$$ \cl{A}=span( \cl{M}^*\cl{B}\cl{M})^{-w^*} \;\;\text{and }\;\;\cl{B}=span( \cl{M}\cl{A}\cl{M}
^*)^{-w^*}.$$
In this case we call the algebras $\cl{A}$ and  $\cl{B}$ {\em TRO equivalent}.
Note that TRO equivalence is a generalisation of unitary equivalence. We show (section 2)
that it is indeed an equivalence relation.

 In the selfadjoint case, TRO equivalence coincides with the existence of an "equivalence bimodule"  
for the algebras (see section 2) thus TRO equivalence implies "Morita equivalence" 
in the sence of Rieffel \cite{rif}. This crucial observation partly motivated our work. 
 In a companion paper we prove that TRO equivalence is the appropriate context in which Rieffel's 
 theory for Morita equivalence of $W^*$-algebras can be generalized to the class of 
possibly nonselfadjoint (abstract) dual operator algebras \cite{ele2}. Using results 
of the present paper this theory is applied in \cite{elejot} to the class 
of reflexive algebras.

Also TRO equivalence is related to the very important notion of stable isomorphism between dual 
operator algebras: In another paper \cite{ep} jointly written with V.I. Paulsen we prove that two 
unital dual operator algebras are stably isomorphic if and only if they have completely isometric normal 
representations with TRO equivalent images. 
These results are generalised, with V.I.Paulsen and I.G. Todorov, to dual operator spaces \cite{ept}.  
 
In the present paper we are concerned with the notion of TRO equivalence within the class
of reflexive (not necessarily selfadjoint) algebras.  
We show (section  3) that  two such algebras are TRO equivalent if and  only  if there exists a
$*$ isomorphism between the commutants of their  diagonals mapping the invariant projection lattice
of the first algebra onto that of the second. This may be thought as a generalisation to the non-selfadjoint case of the remark of Connes \cite{bm} that  two $W^*$
algebras are Morita equivalent in the sense of Rieffel if and  only  if
they have faithful normal representations with isomorphic commutants.

In section  5 we specialise to the case of TRO equivalence of separably acting CSL algebras.
Given the above criterion for TRO equivalence of  reflexive algebras the problem is
the following:

If $\cl{A}, \cl{B}$ are separably acting CSL algebras and
$\phi : \mathrm{Lat}(\cl{A}) \rightarrow \mathrm{Lat}(\cl{B})$
is a lattice isomorphism,
under what conditions does $\phi$ extend to a $*$ isomorphism between the generated von Neumann
algebras  $\mathrm{Lat}(\cl{A})''$ and $\mathrm{Lat}(\cl{B})''$?
The interesting fact is that while $\phi $ always extends to a  $*$ isomorphism between the
generated $C^*$ algebras (Lemma \ref{5.2}), it does not always extend to the $w^*$
closures of these algebras (Remark \ref{4.5}).

In this paper we also consider an equivalence relation
strictly weaker than TRO equivalence, which we call {\em spatial Morita equivalence} (section  4).
Two $w^*$ closed operator algebras $\cl{A}, \cl{B}$ are called {\em spatially Morita equivalent} 
if there exist an $\cl{A}, \cl{B}$ bimodule $\cl{U}$ and a $\cl{B}, \cl{A}$ bimodule $\cl{V}$ such 
that $\cl{A}=span( \cl{U}\cl{V})^{-w^*}$ and  $\cl{B}=span( \cl{V}\cl{U})^{-w^*}$. We show that  
two CSL algebras  
are spatially Morita equivalent if and only if
they have isomorphic lattices. In this case if one of the algebras is ``synthetic" then so
is the other.

\bigskip

We present some definitions and concepts used in this work.

By an algebra $\cl A$ we shall mean an algebra of operators on some Hilbert space;  the \textbf{diagonal} of $\cl A$ is $\Delta (\cl{A})=\cl{A}\cap \cl{A}^*$.

A set of projections of a Hilbert space is called a \textbf{lattice} if it contains the zero and identity projections and is closed under arbitrary suprema and infima. If $\cl{A}$ is a subalgebra of $B(H)$ for some Hilbert space $H$, the set

$$\Lat{\cl A }=\{L\in pr(B(H)): L^\bot \cl{A}L=0\}$$
is a lattice. Dually if $\cl{L}$ is a lattice the space
$$\Alg{\cl L}=\{A\in B(H): L^\bot AL=0 \;\;\forall\;\; L\in \cl{L}\}$$
is an algebra.

A lattice $\cl{L}$ such that $P\in \cl{L}\Leftrightarrow P^\bot \in \cl{L}$ is called an \textbf{ortholattice}. A commutative subspace lattice  \textbf{(CSL)} is  a projection
lattice $\cl{L}$ whose elements commute; the algebra $\mathrm {Alg}(\cl{L})$ is called a 
\textbf{CSL algebra}. A totally ordered CSL is called a \textbf{nest}.

An order-preserving bijection between two lattices is called a \textbf{lattice isomorphism}. If the lattices $\cl{L}_1, \cl{L}_2$ are ortholattices and $\phi :\cl{L}_1\rightarrow \cl{L}_2$ is a lattice isomorphism satisfying $\phi(P^\bot )= \phi(P)^\bot$ for all $P\in \cl{L}_1 $ we call $\phi $ an \textbf{ortholattice isomorphism}.

\medskip

 Let $H_1, H_2$ be Hilbert spaces and $\cl{U}$ a subset of $B(H_1,H_2).$ The \textbf{reflexive hull} of $\cl{U}$ is defined \cite{logshu} to be the space
$$\Ref{ \cl U}=\{T\in B(H_1,H_2): Tx\in \overline{span( \cl{U}x)} \;\text{for\;each}\;x\in H_1\}.$$
Simple arguments show that
\begin{align*}\Ref{\cl U}=&\{T\in B(H_1,H_2)\;\text{for\;all\;projections}\;E,F: E\cl{U}F=0\Rightarrow ETF=0\}\\
=&\{T\in B(H_1,H_2)\;\text{for\;all\;operators}\;A,B: A\cl{U}B=0\Rightarrow ATB=0\}.
\end{align*}
A subspace $\cl{U}$ is called reflexive if $\cl{U}=\Ref{\cl U}.$ A unital algebra is reflexive if and only if $\cl{A}=\mathrm {Alg}(\mathrm {Lat}(\cl{A})).$ CSL algebras are reflexive.
 Every CSL algebra contains a maximal abelian selfadjoint algebra (\textbf{masa} in the sequel). Hence we can view a CSL algebra as a masa bimodule.

Moreover, an algebra is a CSL algebra if and only if it is reflexive and contains a masa.
If $\cl{U}$ is a reflexive masa bimodule, then there exists \cite{dav}, \cite{st} a smallest
 $w^*$ closed masa bimodule which is contained in $\cl{U}$ and whose reflexive hull is the
space $\cl{U}.$ We denote this space by $\cl{U}_{min}.$ Whenever $\cl{U}_{min}=\cl{U}$ we call
the space $\cl{U}$ \textbf{synthetic}. When $\cl{A}$ is a separably acting 
 CSL algebra, the space $\cl{A}_{min}$
 is an algebra which contains the diagonal of $\cl{A}$ and whose lattice is 
$\mathrm {Lat}(\cl{A})$
 \cite{arv}. The first example of a nonsynthetic CSL algebra was given in \cite{arv}.

\medskip

Now we present some concepts introduced in \cite{erd}.

Let $\cl{P}_i=pr(B(H_i)), i=1,2.$ Define $\phi=\mathrm{Map} (\cl{U})$ to be the map $\phi:\cl{P}_1\rightarrow \cl{P}_2$\;which associates to every
$P\in\cl{P}_1$\;the projection onto the subspace
$span( TPy:T\in\;\cl{U},y\in\;H_1)^{-}.$ The map $\phi$\;is $\vee-$continuous (that is, it preserves arbitrary suprema) and $0$ preserving.

Let
$\phi^*=\mathrm{Map} (\mathcal{U}^*),\mathcal{S}_{1,\phi}=\{\phi^*(P)^{\bot}:P\in\cl{P}_2\},\mathcal{S}_{2,\phi}=\{\phi(P):P\in\cl{P}_1\}$ and observe that $\cl{S}_{1,\phi }=\cl{S}_{2,\phi ^*}^\bot .$
Erdos proved that $\mathcal{S}_{1,\phi}$\;is $\wedge $-complete and
contains the identity projection, $\mathcal{S}_{2,\phi}$\;is $\vee $-complete and contains the zero projection, while
$\phi|_{\mathcal{S}_{1,\phi}}:\mathcal{S}_{1,\phi}\rightarrow\mathcal{S}_{2,\phi}$\;is
a bijection. We call the families $\cl{S}_{1,\phi }, \cl{S}_{2,\phi}$ the \textbf{semilattices} of $\cl{U}.$

In fact $$\mathrm{Ref} (\cl{U})=\{T\in\;B(H_1,H_2):\phi(P)^{\bot}TP=0\; \text{for
each}\;P\in\mathcal{S}_{1,\phi}\}.$$
When $\phi(I)=I$ and $\phi^*(I)=I $ we call the space $\cl{U}$ \textbf{essential}.

\medskip

In \cite{kt} it is proved that a TRO $\cl{M}$ is $w^*$ closed if and only if  it is $wot$ closed if and only if  it is
reflexive. In this case, if $\chi =\Map{\cl M}$
$$\cl{M}=\{T\in B(H_1,H_2): TP=\chi (P)T\;\text{for\;all\;}P\in \cl{S}_{1,\chi }\}.$$
 In the following theorem we isolate some consequences of \cite[Theorem 2.10]{kt}.

\begin{theorem}\label{1.1}(i) A TRO $\cl{M}$ is essential if and only if the algebras 
$span( \cl{M}^*\cl{M})^{-w^*},$\\$ span( \cl{M}\cl{M}^*)^{-w^*}$ contain the identity operators.

(ii) If $\cl{M}$ is an essential TRO and $\chi =\Map{\cl M}$ then
$\cl{S}_{1,\chi }=pr((\cl{M}^*\cl{M})^\prime),\\ \cl{S}_{2,\chi }=pr((\cl{M}\cl{M}^*)
^\prime)$ and the map $\chi |_{\cl{S}_{1,\chi }}: \cl{S}_{1,\chi }\rightarrow \cl{S}_
{2,\chi }$ is an ortholattice isomorphism with inverse $\chi^* |_{\cl{S}_{2,\chi }}.$

\end{theorem}

 If $\cl{L}\subset B(H)$ we 
denote by $\cl{L}^\prime$ the set of operators which commute with the elements of 
$\cl{L}$ and the set of projections in $\cl{L}$ by $pr(\cl{L}).$

\section{TRO equivalent algebras}

\begin{definition}\label{2.1.d} Let $\cl{A}, \cl{B}$ be $w^*$ closed algebras 
acting on Hilbert spaces $H_1$ and $H_2$ respectively. If there exists a TRO

$\cl{M}\subset B(H_1,H_2)$ such that $\cl{A}=span( \cl{M}^*\cl{B}\cl{M})^{-w^*}\;\; 
\text{and}\;\;
\cl{B}=span( \cl{M}\cl{A}\cl{M}^*)^{-w^*}$ we write $\cl{A} 
\stackrel{\cl{M}}{\sim}\cl{B}.$ We say that the algebras $\cl{A}, \cl{B}$
 are \textbf{TRO equivalent} if there exists a TRO 
$\cl{M}$ such that $\cl{A} 
\stackrel{\cl{M}}{\sim}\cl{B}.$ 
\end{definition}

A simple example of TRO equivalent, not necessarily selfadjoint 
algebras, is the following. Take a unital $w^*$ closed  algebra $\cl A\subset B(H)$ 
and let $$\cl{B}= \left[\begin{array}{clr}\cl{A} & \cl{A}\\ \cl{A} & \cl{A}\end{array}
\right] \subset B(H\oplus H),\;\; \cl{M}=\left[\begin{array}{clr}\Delta (\cl{A}) 
\\ \Delta (\cl{A})\end{array}\right]\subset B(H, H\oplus H).$$ It is easy to see 
that $\cl{A} \stackrel{\cl{M}}{\sim}\cl{B}.$

\begin{proposition}\label{2.1} Let $\cl{A}\subset B(H_1), \cl{B}\subset B(H_2)$ 
be $w^*$ closed algebras. The following are equivalent:

 (i) The algebras $\cl{A}, \cl{B}$ are TRO equivalent.

(ii) There exists an essential TRO $\cl{M}\subset B(H_1,H_2)$ such that
$\cl{M}^*\cl{B}\cl{M}\subset \cl{A}$ and $\cl{M}\cl{A}\cl{M}^*\subset \cl{B}.$

 If (ii) holds then $\cl{A} \stackrel{\cl{M}}{\sim}\cl{B}.$
\end{proposition}
\textbf{Proof} Let $\cl{N}\subset B(H_1,H_2)$ be a TRO such that 
$\cl{A}=span( \cl{N}^*\cl{B}\cl{N})^{-w^*}$\\ $ \text{and}\;\;
\cl{B}=span( \cl{N}\cl{A}\cl{N}^*)^{-w^*}.$

If $P$ is the projection onto $\ker \cl N$ and $Q$ is the projection onto $\ker \cl N^*$, it
is clear that $\cl M \equiv \cl N + QB(H_1,H_2)P$ is an essential TRO such that
$\cl{M}^*\cl{BM}\subset \cl{A}$ and $\cl{MAM}^*\subset \cl{B}$.

Conversely, if $\cl{M}$ is an essential TRO such that $\cl{M}^*\cl{B}\cl{M}\subset 
\cl{A}$ and $\cl{M}\cl{A}\cl{M}^* \subset \cl{B}$, then $span( \cl{M}^*\cl{B}\cl{M})^{-w^*}
\subset \cl{A}$ and, since
$\cl{M}\cl{A}\cl{M}^*\subset \cl{B}$, we have
\begin{align*}&\cl{M}^*\cl{M}\cl{A}\cl{M}^*\cl{M}\subset \cl{M}^*\cl{B}\cl{M}\Rightarrow \\ 
& span( \cl{M}^*\cl{M})^{-w^*}\cl{A}span( \cl{M}^*\cl{M})^{-w^*}\subset span( \cl{M}^*\cl{B}
\cl{M})^{-w^*}.\end{align*}
Since $I\in span( \cl{M}^*\cl{M})^{-w^*}$ it follows that
$\cl{A}\subset span( \cl{M}^*\cl{B}\cl{M})^{-w^*}.$

 Symmetrically we have $\cl{B}=span( \cl{M}\cl{A}\cl{M}^*)^{-w^*}. \qquad  \Box$

\bigskip

 We wish to prove that TRO equivalence is an equivalence relation. We need the following lemma.

\begin{lemma}\label{2.2} Let $\cl{S}_i$ be a set of projections on the Hilbert space 
$H_i, i=1,2$, $\chi : \cl{S}_1\rightarrow \cl{S}_2$ a map onto $\cl{S}_2,$ and 
$$\cl{M}=\{T\in B(H_1,H_2): TL=\chi (L)T\;\; \text{for\;all\;}L\in \cl{S}_1\}.$$ 
Observe that the space $\cl{M}$ is a reflexive TRO \cite{kt}. Moreover, if we have 
the information that it is essential, then 
$$span( \cl{M}^*\cl{M})^{-w^*}=(\cl{S}_1)^\prime\;\;\text{and}\;\; 
span( \cl{M}\cl{M}^*)^{-w^*}=(\cl{S}_2)^\prime.$$
\end{lemma}
\textbf{Proof}\;\;Let $\phi =\mathrm {Map}(\cl{M}).$ We can observe that $\cl{M}(\cl{S}_1)^\prime\subset \cl{M},$ so if $P$ is a projection then $(\cl{S}_1)^\prime\cl{M}^*P(H_2)\subset \cl{M}^*P(H_2)$ therefore $(\cl{S}_1)^\prime\phi ^*(P)(H_1)\subset \phi ^*(P)(H_1).$ Since $\cl S_1$ is selfadjoint, it follows that $\phi ^*(P)\in (\cl{S}_1)^{\prime\prime}.$

We proved that $\cl{S}_{2,\phi ^*}\subset (\cl{S}_1)''$; thus $\cl{S}_{1,\phi}\subset 
(\cl{S}_1)''$ and so $(\cl{S}_{1,\phi})^\prime\supset  (\cl{S}_1)^{\prime}.$ But $
span( \cl{M}^*\cl{M})^{-w^*}=(\cl{S}_{1,\phi} )^\prime$ (Theorem \ref{1.1}) since 
$\cl{M}$ is an essential TRO. We proved that $
span( \cl{M}^*\cl{M})^{-w^*}\supset (\cl{S}_1)^\prime.$

Clearly, $\cl{M}^*\cl{M}\subset (\cl{S}_1)'$ and so $span( \cl{M}^*\cl{M})^{-w^*}=(\cl{S}_1)^\prime.$

Since $\chi$ maps onto $\cl S_2$ we see that $\cl{M}^*(\cl{S}_2)'\subset \cl M^*$ and
similar arguments show that $span( \cl{M}\cl{M}^*)^{-w^*}= (\cl{S}_2)^\prime.
\qquad \Box$

\begin{theorem}\label{2.3} TRO equivalence is an equivalence relation.
\end{theorem}
\textbf{Proof} We only have to prove  transitivity. Let $\cl{A}, \cl{B}, \cl{C}$ be $w^*$ 
closed algebras, acting on the Hilbert spaces $H_1, H_2, H_3$ respecively, and 
$\cl{M}, \cl{N}$ be essential TRO's such that $\cl{B} \stackrel{\cl{M}}{\sim}\cl{A}$ and 
$\cl{B} \stackrel{\cl{N}}{\sim}\cl{C}.$ Thus
\begin{align*}
span( \cl{MBM}^*)^{-w^*}=\cl A,& \;\; span( \cl{NBN}^*)^{-w^*}=\cl C \\
\text{and } \;\;\;
span( \cl M^*\cl{AM})^{-w^*}=\cl B &=span( \cl N^*\cl{CN})^{-w^*}.
\end{align*}
Define $$\cl S=pr((\cl{M}^*\cl{M})'\cap (\cl{N}^*\cl{N})')$$
and note that $$\cl S'=((\cl{M}^*\cl{M})\cup (\cl{N}^*\cl{N}))''.$$ 
Since $\cl{M}^*\cl{M}\cl{B}
\cl{M}^*\cl{M}\subset \cl B$ and similarly for $\cl{N}$ it follows that 
$\cl{S}^\prime\cl{B}\cl{S}^\prime\subset \cl{B}.$
Let $\chi=\Map{ \cl M}$ and $\phi=\Map{ \cl N}$. 
Define the TRO's
\begin{align*}
& \cl Z = \{T\in B(H_2, H_1): TL=\chi(L)T, \, L\in \cl S\}\\ &
\cl Y = \{T\in B(H_2, H_3): TL=\phi(L)T, \, L\in \cl S\}.\end{align*}
The map $\chi$ is an ortholattice isomorphism from $pr((\cl{M}^*\cl{M})^\prime)$ onto 
$pr((\cl{M}\cl{M}^*)^\prime)$ (Theorem \ref{1.1}). 
Since $\cl S\subset pr((\cl{M}^*\cl{M})')$ it follows that $\cl M\subset \cl Z$.
Similarly $\cl N\subset \cl Y$ and thus both $\cl Z$ and $\cl Y$ are essential TRO's.
From the previous lemma we obtain
$$span( \cl{Y}^*\cl{Y})^{-w^*}=\cl{S}'=span( \cl{Z}^*\cl{Z})^{-w^*}.$$

We claim that $$span( \cl Z^*\cl{AZ})^{-w^*}=\cl B \;\;
\text{and }\;\; span( \cl{ZBZ}^*)^{-w^*}=\cl{A}.$$
Indeed since $\cl{Z}^*\cl{Z}\subset \cl{S}'$ and $\cl M\subset \cl Z$ we have
\begin{align*} \cl{Z}^*\cl{ZBZ}^*\cl{Z}\subset \cl{B} \Rightarrow & \cl{M}\cl{Z}^*(\cl{ZBZ}^*)\cl{Z}\cl{M}^*\subset
\cl{MBM}^*\subset \cl{A} \\
 \Rightarrow & \cl{MM}^*(\cl{ZBZ}^*)\cl{MM}^*\subset \cl{A}. \end{align*}
Since $\cl M$ is essential and $span(\cl{MM}^*)$ (resp. $span(\cl{M}^*\cl M)$) is a $*-$algebra, 
one can find a bounded net in $span(\cl{MM}^*)$ (resp. $span(\cl{M}^*\cl M)$)  converging 
strongly to the identity operator. Since $\cl A$ is $w^*-$closed it follows that 
 $\cl{ZBZ}^* \subset \cl{A}$ and hence  $span(\cl{ZBZ}^*)^{-w^*} \subset \cl{A}.$

On the other hand $\cl{A}=span( \cl{MBM^*})^{-w^*}\subset span( \cl{ZBZ^*})^{-w^*}$ hence
$\cl{A}=span( \cl{ZBZ}^*)^{-w^*}.$ Therefore
$\cl{Z}^*\cl{AZ}=\cl{Z}^*span( \cl{ZBZ}^*)^{-w^*}\cl{Z}\subset
\cl{B}$ because $\cl{Z}^*\cl{ZBZ}^*\cl{Z}\subset \cl{B}.$
It follows by Proposition \ref{2.1} that $\cl{B}=span( \cl{Z}^*\cl{AZ})^{-w^*}.$

In the same way we have
$$span( \cl Y^*\cl{CY})^{-w^*}=\cl B \;\;
\text{and }\;\; span( \cl{YBY}^*)^{-w^*}=\cl{C}.$$

Now put  $\cl{L}=span( \cl{YZ}^*)^{-w^*}.$ Since
$span( \cl{Y}^*\cl{Y})^{-w^*}=\cl{S}'=span( \cl{Z}^*\cl Z)^{-w^*}$
we have $\cl{YZ}^*\cl{Z}\subset \cl{Y}.$ It follows that
$$(\cl{YZ}^*)(\cl{YZ}^*)^*(\cl{YZ}^*)=\cl{YZ}^*\cl{ZY}^*\cl{YZ}^*\subset
\cl{YY}^*\cl{YZ}^* \subset \cl{YZ}^*$$ since $\cl Y$ is a TRO, hence
$\cl{LL}^*\cl{L}\subset \cl{L}.$ Thus the space
$\cl{L}$ is a TRO and it is essential because the spaces
$\cl{Y}$ and $\cl{Z}^*$ are essential TRO's.

To complete the proof it remains to show that
$$\cl L^*\cl{CL}\subset \cl A \;\;
\text{and }\;\; \cl{LAL}^*\subset\cl{C}.$$
Indeed, since $\cl Y^*\cl{CY}\subset \cl B$ we have
$\cl{ZY}^*\cl{CYZ}^*\subset \cl{ZBZ}^*\subset \cl A$
and thus $\cl{L}^*\cl{CL}\subset \cl{A}$.
On the other hand,  $\cl{YZ}^*\cl{AZY}^*\subset
\cl{YBY}^*\subset\cl{C}$ and therefore
$\cl{LAL}^*\subset \cl C$.  $\qquad  \Box$

\begin{remark}\label{2.4}From the previous proof it follows that, if $\cl{A}, \cl{B}$
are TRO equivalent algebras and $\cl{B}, \cl{C}$ are also TRO equivalent algebras, then
there exist  essential TRO's $\cl{Z}, \cl{Y}$ generating the same von Neumann algebra
such that $\cl{B} \stackrel{\cl{Z}}{\sim}\cl{A},$ $\cl{B} \stackrel{\cl{Y}}{\sim}\cl{C}$ 
and the space $\cl{L}=span( \cl{YZ}^*)^{-w^*}$ is an essential TRO satisfying
 $\cl{A} \stackrel{\cl{L}}{\sim}\cl{C}.$
\end{remark}

\begin{proposition}\label{2.5}Let $\cl{A}, \cl{B}$ be $w^*$ closed algebras and $\cl{M}$ an essential TRO such that $\cl{A} \stackrel{\cl{M}}{\sim}\cl{B}.$ Then $\Delta (\cl{A}) \stackrel{\cl{M}}{\sim}\Delta (\cl{B}).$
\end{proposition}
\textbf{Proof} It is obvious that $\cl{M}^*\Delta (\cl{B})\cl{M}\subset \Delta (\cl{A})$ and $\cl{M}\Delta (\cl{A})\cl{M}^*\subset \Delta (\cl{B}).$ By Proposition \ref{2.1} follows that $\Delta (\cl{A}) \stackrel{\cl{M}}{\sim}\Delta (\cl{B}). \qquad  \Box$

\begin{lemma}\label{2.6} Let $\cl{A}, \cl{B}$ be unital $w^*$ closed algebras, 
$\cl{M}$ be an essential TRO such that $\cl{A} \stackrel{\cl{M}}{\sim}\cl{B}$ 
and $\chi =\mathrm {Map}(\cl{M}).$ Then $\mathrm{Ref}(\cl{A}) \stackrel{\cl{M}}{\sim}\mathrm{Ref}
(\cl{B}).$ Also the map $\chi : 
pr(\Delta (\cl{A})^\prime)\rightarrow pr(\Delta (\cl{B})^\prime)$ is an 
orthoisomorphism and $\chi (\mathrm {Lat}(\cl{A}))=\mathrm {Lat}(\cl{B}).$
\end{lemma}
\textbf{Proof} By the above proposition $\Delta (\cl{A}) \stackrel{\cl{M}}{\sim}\Delta 
(\cl{B}).$ From \cite[Corollary 5.9]{kt} it follows that $\chi(pr(\Delta (\cl{A})^\prime))=
pr(\Delta (\cl{B})^\prime).$ Since $\cl{M}^*\cl{M}\subset \Delta (\cl{A})$ we have 
$pr ((\cl{M}^*\cl{M})^\prime)\supset pr ((\Delta (\cl{A}))^\prime).$ So by Theorem 
\ref{1.1} the map $\chi : pr(\Delta (\cl{A})^\prime)\rightarrow pr(\Delta (\cl{B})^\prime)$ 
is an orthoisomorphism.

 If $E,F$ are projections such that $E\cl{A}F=0$ then $E\cl{M}^*\cl{B}\cl{M}F=0$ so 
$E\cl{M}^*\mathrm{Ref}(\cl{B})\cl{M}F=0.$ It follows that $$\cl{M}^*\mathrm{Ref}
(\cl{B})\cl{M}\subset \mathrm{Ref}(\cl{A}).$$ Similarly we can prove that 
$\cl{M}\mathrm{Ref}(\cl{A})\cl{M}^*\subset \mathrm{Ref}(\cl{B}),$ hence 
$\mathrm{Ref}(\cl{A}) \stackrel{\cl{M}}{\sim}\mathrm{Ref}(\cl{B}).$ Since 
$\mathrm{Lat}(\mathrm{Ref}(\cl{A}))=\mathrm{Lat}(\cl{A})$ and similarly for 
$\cl{B},$ using again \cite[Corollary 5.9]{kt} we have $\chi (\mathrm {Lat}
(\cl{A}))=\mathrm {Lat}(\cl{B}). \qquad  \Box$

\begin{remark}\label{2.7}\em{(i) Let $\cl{A}, \cl{B}$ be TRO equivalent 
$w^*$-closed algebras and suppose that the algebra $\cl{A}$ is reflexive. 
Then the algebra $\cl{B}$ is reflexive. Indeed if $\cl{M}$ is an essential 
TRO such that $\cl{A} \stackrel{\cl{M}}{\sim}\cl{B}$ as in the proof of 
Lemma \ref{2.6} it follows that $\cl{A} \stackrel{\cl{M}}{\sim}\mathrm{Ref}
(\cl{B}),$ hence $\Ref{\cl{B}}=\cl{B}.$

(ii) An orthoisomorphism $\chi: pr(\cl{C})\rightarrow pr(\cl{D})$, where 
$\cl{C}$ and $\cl{D}$ are von Neumann algebras, does not necessarily extend 
to a $*-$homomorphism between the algebras. For example choose \cite{kr} 
nonabelian $*$ anti-isomorphic von Neumann algebras $\cl{C}, \cl{D},$ 
$\theta : \cl{C}\rightarrow \cl{D}$ a $*$ anti-isomorphism  and let 
$\chi =\theta |_{pr(\cl{C})}.$ Compare now Lemma \ref{2.6} and Theorem \ref{3.3}.}
\end{remark}

\begin{proposition}\label{2.8} Let $\cl{A}, \cl{B}$ be unital $w^*$ closed algebras 
acting on the Hilbert spaces $H_1, H_2$ respectively and $\cl{M}$ be an essential TRO such that $\cl{A} \stackrel{\cl{M}}{\sim}\cl{B}.$ Then there exists a TRO $\cl{N}$ which contains  $\cl{M}$ such that $\cl{A} \stackrel{\cl{N}}{\sim}\cl{B}$ and $\Delta 
(\cl{A})=span( \cl{N}^*\cl{N})^{-w^*}, \Delta (\cl{B})=span( \cl{N}\cl{N}^*)^{-w^*}.$
\end{proposition}
\textbf{Proof} Let $\chi =\mathrm {Map}(\cl{M}).$ From Lemma \ref{2.6} follows that $\chi (pr(\Delta (\cl{A})^\prime)=pr(\Delta (\cl{B})^\prime).$ Let
$$\cl{N}=\{T\in B(H_1, H_2): TL=\chi (L)T\;\; \text{for\;all\;}L\in pr(\Delta (\cl{A})^\prime)\}.$$
Since $\cl{S}_{1,\chi }=pr((\cl{M}^*\cl{M})^\prime)\supset
pr(\Delta (\cl{A})^\prime)$ we have that $\cl{M}\subset \cl{N}$ so
the TRO $\cl{N}$ is essential.

Using Lemma \ref{2.2} we have $\Delta (\cl{A})=span( \cl{N}^*\cl{N})^{-w^*}, \Delta (\cl{B})=span( \cl{N}\cl{N}^*)^{-w^*}.$ We shall show that $\cl{A}=span( \cl{N}^*\cl{B}\cl{N})^{-w^*}.$ Since $\cl{M}^*\cl{B}\cl{M}\subset \cl{N}^*\cl{B}\cl{N},$ we get $\cl{A}\subset span( \cl{N}^*\cl{B}\cl{N})^{-w^*}.$ Now we have \begin{align*}\cl{N}\cl{N}^*\cl{B}\cl{N}\cl{N}^*\subset \cl{B}\Rightarrow & \cl{M}^*\cl{N}\cl{N}^*\cl{B}\cl{N}\cl{N}^*\cl{M}\subset \cl{M}^*\cl{B}\cl{M}\subset \cl{A} \\ \Rightarrow & \cl{M}^*\cl{M}\cl{N}^*\cl{B}\cl{N}\cl{M}^*\cl{M}\subset \cl{A} \end{align*} hence $span( \cl{M}^*\cl{M})^{-w^*}\cl{N}^*\cl{B}\cl{N}span( \cl{M}^*\cl{M})^{-w^*}\subset \cl{A}.$ \\
It follows that $\cl{N}^*\cl{B}\cl{N}\subset \cl{A}$ therefore $\cl{A}=span( \cl{N}^*\cl{B}\cl{N})^{-w^*}.$ Now since $\cl{N}\cl{N}^*\subset \cl{B}$ we have $\cl{N}\cl{A}\cl{N}^*=\cl{N}span( \cl{N}^*\cl{B}\cl{N})^{-w^*}\cl{N}^*\subset \cl{B}.$ We proved that $\cl{A} \stackrel{\cl{N}}{\sim}\cl{B}.\qquad \Box$

\bigskip

We isolate the following consequence of the above proposition:

\begin{corollary} If the unital $w^*$ closed algebras $\cl{A}, \cl{B}$ are TRO equivalent 
then the diagonal algebras $\Delta (\cl{A}), \Delta (\cl{B})$ are Morita equivalent 
in the sense of Rieffel \cite{rif}.
\end{corollary}

The following proposition says that if two non-unital
$w^*$ closed algebras are TRO equivalent, there exist TRO
equivalent unital algebras which contain the previous algebras as
ideals.

\begin{proposition}\label{2.9} Let $\cl{A}, \cl{B}$ be non unital $w^*$ closed algebras and $\cl{M}$ be an essential TRO such that $\cl{A} \stackrel{\cl{M}}{\sim}\cl{B}.$ If $\cl{A}_{\cl{M}}=span( \cl{A},\cl{M}^*\cl{M})^{-w^*}, \cl{B}_{\cl{M}}=span( \cl{B},\cl{M}\cl{M}^*)^{-w^*},$ then

(i)The spaces $\cl{A}_{\cl{M}},\cl{B}_{\cl{M}}$ are unital algebras.

(ii)The algebra $\cl{A}$ (respectively $\cl{B}$) is an ideal of $\cl{A}_{\cl{M}}$ (respectively $\cl{B}_{\cl{M}}$).

(iii)$\cl{A}_{\cl{M}} \stackrel{\cl{M}}{\sim}\cl{B}_{\cl{M}}.$

(iv)There exists an essential TRO $\cl{N}$ which contains $\cl{M}$ such that 
$\cl{A} \stackrel{\cl{N}}{\sim}\cl{B},$ $\Delta
(\cl{A}_{\cl{M}})=span( \cl{N}^*\cl{N})^{-w^*}, \Delta
(\cl{B}_{\cl{M}})=span( \cl{N}\cl{N}^*)^{-w^*}.$ (Observe that
$\cl{A}_{\cl{M}}=\cl{A}_{\cl{N}}$ and
$\cl{B}_{\cl{M}}=\cl{B}_{\cl{N}}).$
\end{proposition}
\textbf{Proof} Claims (i),(ii) are consequences of the relations $\cl{A}\cl{M}^*\cl{M}\subset \cl{A},$\linebreak $\cl{M}^*\cl{M}\cl{A}\subset \cl{A},\;\; \cl{B}\cl{M}\cl{M}^*\subset \cl{B},\;\; \cl{M}\cl{M}^*\cl{B}\subset \cl{B}.$\\
Also, since $\cl{M}\cl{A}\cl{M}^*\subset \cl{B}$ and $\cl{M}(\cl{M}^*\cl{M})\cl{M}^*\subset \cl{M}\cl{M}^*$ it easily follows that $\cl{M}\cl{A}_{\cl{M}}\cl{M}^*\subset \cl{B}_{\cl{M}}.$ Similarly we get $\cl{M}\cl{B}_{\cl{M}}\cl{M}^*\subset \cl{A}_{\cl{M}}.$

(iv) Since $\cl{A}_{\cl{M}} \stackrel{\cl{M}}{\sim}\cl{B}_{\cl{M}},$ by the 
previous proposition there exists an essential TRO $\cl{N}$ containing $\cl{M}$ 
such that $\cl{A}_{\cl{M}} \stackrel{\cl{N}}{\sim}\cl{B}_{\cl{M}}$ and 
$$\Delta (\cl{A}_{\cl{M}})=span( \cl{N}^*\cl{N})^{-w^*}, \Delta 
(\cl{B}_{\cl{M}})=span( \cl{N}\cl{N}^*)^{-w^*}.$$
It remains to show that $\cl{A} \stackrel{\cl{N}}{\sim}\cl{B}.$\\
Since $\cl{N}\cl{N}^*\subset \cl{B}_{\cl{M}}$ and  $\cl{B}$ is an
ideal of $\cl{B}_{\cl{M}}$ we have  \begin{align*} &\cl{N}\cl{N}^*\cl{B}
\cl{N}\cl{N}^*\subset \cl{B} \Rightarrow \cl{M}^*\cl{N}\cl{N}^*\cl{B}
\cl{N}\cl{N}^*\cl{M}\subset \cl{M}^*\cl{B}\cl{M}\subset \cl{A} \\ \Rightarrow & \cl{M}^*\cl{M}\cl{N}^*\cl{B} \cl{N}\cl{M}^*\cl{M} \subset
\cl{A} \Rightarrow span( \cl{M}^*\cl{M})^{-w^*}\cl{N}^*\cl{B}
\cl{N}span( \cl{M}^*\cl{M})^{-w^*}\subset \cl{A} \end{align*} hence
$\cl{N}^*\cl{B} \cl{N}\subset \cl{A}.$ Similarly we can prove
that $\cl{N}^*\cl{A} \cl{N}\subset \cl{B}. \qquad \Box$

\begin{proposition}\label{2.10} Let $\cl{A}, \cl{B}$ be $w^*$ closed algebras and 
$\cl{M}$ be an essential TRO such that $\cl{A} \stackrel{\cl{M}}{\sim}\cl{B}.$ 
If $\cl{J}$ is a $w^*$ closed $\cl{A}-$bimodule then the space $F(\cl{J})=
span( \cl{M}\cl{J}\cl{M}^*)^{-w^*}$ is a $\cl{B}-$bimodule and  $\cl{J} \stackrel
{\cl{M}}{\sim}F(\cl{J}).$ The map $F$ is a bijection between $w^*$ closed 
bimodules of $\cl{A}$ and those of $\cl{B}.$ Moreover the restriction of $F$ 
to the set of two sided $w^*$ closed ideals of $\cl{A}$ maps onto those of $\cl{B}.$
\end{proposition}
\textbf{Proof} Since $\cl{A}\cl{M}^*\cl{M}\subset \cl{A}$ we have 
$\cl{M}\cl{A}\cl{M}^*\cl{M}\cl{J}\cl{M}^*\subset \cl{M}\cl{J}\cl{M}^*.$ 
Hence $\cl{B}F(\cl{J})\subset F(\cl{J}).$ Similarly we have $F(\cl{J})\cl{B}
\subset F(\cl{J}).$ \\
If $\cl{I}$ is a $w^*$ closed bimodule of $\cl{B},$ the space $\cl{J}=
span( \cl{M}^*\cl{I}\cl{M})^{-w^*}$ is a bimodule of $\cl{A}$ and 
$F(\cl{J})=\cl{I}.$ So the map $F$ is onto. Clearly, F is an injection. 
Also observe that if $\cl{J}$ is a two sided $w^*$ closed ideal of $\cl{A}$ 
then the space $F(\cl{J})$ is a two sided $w^*$ closed ideal of $\cl{B}.$ $\qquad \Box$

\medskip

The following proposition is proved easily.

\begin{proposition}\label{2.11} Let $\cl{A}, \cl{B}$ be $w^*$ closed algebras and 
$\cl{M}$ be an essential TRO such that $\cl{A} \stackrel{\cl{M}}{\sim}\cl{B}.$ 
We denote by $\cl{K}(\cl{A})$ (respectively $\cl{K}(\cl{B})$) the set of compact 
operators in $\cl{A}$ (resp. $\cl{B}$), by $\cl{F}(\cl{A})$ (resp. $\cl{F}(\cl{B})$) 
the set of finite rank operators in $\cl{A}$ (resp. $\cl{B}$), by $R_1(\cl{A})$ 
(resp. $R_1(\cl{B})$) the set of rank 1 operators in $\cl{A}$(resp. $\cl{B}$). 
Then it follows \begin{align*}& \cl{K}(\cl{A})^{-w^*}\stackrel{\cl{M}}{\sim}\cl{K}(\cl{B})^{-w^*},
\;\; \cl{F}(\cl{A})^{-w^*}\stackrel{\cl{M}}{\sim}\cl{F}(\cl{B})^{-w^*},\\  
& span( R_1(\cl{A}))^{-w^*}\stackrel{\cl{M}}{\sim}span( R_1(\cl{B}))^{-w^*}.\end{align*}
\end{proposition}

\section{TRO equivalent reflexive algebras}

The goal of this section is to determine sufficient and necessary conditions for TRO equivalence of reflexive algebras. The following lemma is known. See for example \cite[8.5.32]{bm}. We include a proof for completeness.

\begin{lemma}\label{3.1} Let $\cl{C}, \cl{E}$ be von Neumann algebras acting on the Hilbert spaces $H_1, H_2$ respectively 
, $\theta :\cl{C}\rightarrow \cl{E}$ be a $*$ isomorphism and 
$$\cl{M}=\{T\in B(H_1, H_2): 
TA=\theta (A)T\;\; \text{for\;all\;}A\in \cl{C}\}.$$ Then the space $\cl{M}$ is an essential TRO.
\end{lemma}
\textbf{Proof} Let $\cl{D}=\{A\oplus \theta (A): A\in \cl{C}\}.$ Since $\theta$ is $w^*$ continuous, as a $*$ isomorphism between von Neumann algebras \cite[I.4.3, Corollaire 2]{dix}, the space $\cl{D}$ is a von Neumann algebra. The commutant of $\cl{D}$ is the algebra $$\left[\begin{array}{clr}\cl{C}^\prime&\cl{M}^*\\\cl{M}&\cl{E}^\prime\end{array}\right].$$
Let $\phi =\mathrm {Map}(\cl{M}).$ Since $\cl{E}^\prime\cl{M}\subset \cl{M}$ we have that 
$\phi (I)^\perp \cl{E}^\prime\phi (I)=0,$ hence $\phi (I)\in \cl{E}.$ 
Let $P=0\oplus \phi (I)^\perp .$ Since $\phi (I)\cl{M}=\cl{M}$ and $\phi (I)\in \cl{E}$ we can verify that $P^\perp \cl{D}^\prime P=0,$ hence $P\in \cl{D}.$ It follows that the projection $P$ is of the form $A \oplus \theta (A)$ for an operator $A \in \cl{C}.$ Thus $\phi (I)=I.$ Similarly we can prove that $\phi ^*(I)=I,$ so the space $\cl{M}$ is an essential TRO.$\qquad \Box$

\medskip

We give a new proof of Connes' remark (see the introduction) and also show that the
isomorphism between the commutants extends the map of the Morita equivalence bimodule.
This fact will be useful below.

\begin{theorem}\label{3.2} Let $\cl{A}, \cl{B}$ be von Neumann algebras acting on 
the Hilbert spaces $H_1, H_2$ respectively , $\cl{M}$ be an essential TRO such that 
$\cl{A}=span( \cl{M}^*\cl{M})^{-w^*},$\\$ \cl{B}=span( \cl{M}\cl{M}^*)^{-w^*}$ and $\chi =\mathrm {Map}(\cl{M}).$ Then there exists a $*$ isomorphism $\theta : \cl{A}^\prime\rightarrow \cl{B}^\prime$ which extends the map $\chi |_{pr(\cl{A}^\prime)}.$\\ Conversely if the algebras $\cl{A}^\prime, \cl{B}^\prime$ are $*$ isomorphic, the algebras $\cl{A}, \cl{B}$ are TRO equivalent.
\end{theorem}
\textbf{Proof}  By Theorem \ref{1.1},
$$\cl{M}=\{T\in B(H_1, H_2): TL=\chi (L)T\;\; \text{for\;all\;}L\in pr(\cl{A}^\prime)\}$$
Let $\cl{L}=\{L\oplus \chi (L): L\in pr(\cl{A}^\prime)\}.$ We can verify that $$\cl{C}=\cl{L}^\prime=\left[\begin{array}{clr}\cl{A}&\cl{M}^*\\\cl{M}&\cl{B}\end{array}\right].$$
So the algebra $\cl{C}$ is a von Neumann algebra acting on the direct sum of the corresponding Hilbert spaces.

An easy calculation shows that the commutant of $\cl{C}$ is the set 
$$\left\{\left[\begin{array}{clr}T&0\\0&S\end{array}\right]: T\in \cl{A}^\prime, 
S\in \cl{B}^\prime\;\;\text{such\;that\;} SM=MT\;\forall M\in \cl{M} \right\}.$$
Let $$\pi _1: \cl{C}^\prime\rightarrow \cl{A}^\prime: \left[\begin{array}
{clr}T&0\\0&S\end{array}\right]\rightarrow T $$
  $$\pi _2: \cl{C}^\prime\rightarrow \cl{B}^\prime: \left[\begin{array}{clr}T&0\\0&S\end
{array}\right]\rightarrow S. $$ We shall show that the maps $\pi_1, \pi_2 $ are surjective.
 Clearly $\pi _1(\cl{C}^\prime)$ is a von Neumann algebra, so it suffices to show that 
$\pi _1(\cl{C}^\prime)^\prime\subset \cl{A}.$\\
If $A\in \pi _1(\cl{C}^\prime)^\prime$ then $AT=TA$ for all $\left[\begin{array}{clr}
T&0\\0&S\end{array}\right]\in \cl{C}^\prime.$ Thus $$\left[\begin{array}{clr}A&0\\0&0\end{array}
\right]\left[\begin{array}{clr}T&0\\0&S\end{array}\right]=\left[\begin{array}{clr}T&0\\0
&S\end{array}\right]\left[\begin{array}{clr}A&0\\0&0\end{array}\right]\;\;\text{for\; 
all}\;\; \left[\begin{array}{clr}T&0\\0&S\end{array}\right]\in \cl{C}^\prime,$$  
hence $\left[\begin{array}{clr}A&0\\0&0\end{array}\right]\in \cl{C},$ and so $A\in \cl{A}.$

\medskip

If $\left[\begin{array}{clr}T&0\\0&S_1\end{array}\right], \left[\begin{array}{clr}T&0\\0&S_2
\end{array}\right]\in \cl{C}^\prime$ then $S_1M=MT=S_2M$ for all $M\in \cl{M}.$ Since the TRO $\cl{M}$ is essential we have $S_1=S_2.$\\
The conclusion is that we can define a map $\theta : \cl{A}^\prime\rightarrow \cl{B}^\prime$ 
such that \\$\theta (T)=S\Leftrightarrow \left[\begin{array}{clr}T&0\\0&S\end{array}\right]
\in \cl{C}^\prime.$ The map $\theta $ is a $*$ isomorphism. We shall show that $\theta $ 
is an extension of $\chi |_{pr(\cl{A}^\prime)}.$\\
If $L\in pr(\cl{A}^\prime)$ we have $ML=\theta (L)M$ for all $M\in \cl{M}.$ It follows $\theta (L)^\perp \cl{M}L=0$ and from this $\chi (L)\leq \theta (L).$ Also $\theta (L) \cl{M}L^\perp =0$ hence $\chi (L^\perp )\leq \theta (L)^\perp .$ By Theorem \ref{1.1} $\chi(L^\bot )= \chi(L)^\bot $ so $\chi (L)= \theta (L).$

\medskip

Conversely, let $\theta :\cl{A}^\prime\rightarrow \cl{B}^\prime$ be a $*$ isomorphism
and $$\cl{M}=\{T: TA=\theta (A)T\;\; \text{for\;all\;}A\in \cl{A}^\prime\}.$$ The space $\cl{M}$ is an essential TRO by the previous lemma. It is obvious that $\cl{M}^*\cl{B}\cl{M}\subset \cl{A}$ and $\cl{M}\cl{A}\cl{M}^*\subset \cl{B}. \qquad \Box$

\begin{theorem}\label{3.3} Two unital reflexive algebras $\cl{A},\cl{B}$ are TRO equivalent
if and only if there exists a $*$ isomorphism $\theta : \Delta (\cl{A})^\prime
 \rightarrow \Delta (\cl{B})^\prime $ such that $\theta (\mathrm {Lat}(\cl{A}))=\mathrm
 {Lat}(\cl{B}).$
\end{theorem}
\textbf{Proof} Let $\cl{A}, \cl{B}$ be TRO equivalent algebras, 
acting on the Hilbert spaces $H_1, H_2$ respectively . By Proposition \ref{2.8} there exists
 an essential TRO $\cl{M}$ such that
$\cl{A} \stackrel{\cl{M}}{\sim}\cl{B}$ and $\Delta (\cl{A})=span( \cl{M}^*\cl{M})^{-w^*}, \Delta (\cl{B})=span( \cl{M}\cl{M}^*)^{-w^*}.$\\
From the previous theorem there exists a $*$ isomorphism $\theta :\Delta (\cl{A})^\prime\rightarrow \Delta (\cl{B})^\prime$ which extends the map $\chi |_{pr(\Delta (\cl{A})^\prime)}.$ From Lemma \ref{2.6}  $\theta(\mathrm {Lat}(\cl{A}))= \chi(\mathrm {Lat}(\cl{A}))=\mathrm {Lat}(\cl{B}).$\\
Conversely, let $\theta :\Delta (\cl{A})^\prime\rightarrow \Delta (\cl{B})^\prime$ be a $*$ isomorphism such that  $\theta(\mathrm {Lat}(\cl{A}))=\mathrm {Lat}(\cl{B})$ and define
$$\cl{M}=\{T\in B(H_1, H_2): TA=\theta (A)T\;\; \text{for\;all\;}A\in \Delta (\cl{A})^\prime\}.$$
By Lemma \ref{3.1} the space $\cl{M}$ is an essential TRO. It remains to show that $\cl{A} \stackrel{\cl{M}}{\sim}\cl{B}.$ Let $A\in \cl{A}, L\in \mathrm {Lat}(\cl{A}), M_1,M_2\in \cl{M}$ then $M_1AM_2^*\theta (L)=M_1ALM_2^*=M_1LALM_2^*=\theta (L)M_1AM_2^*\theta (L).$ Hence $M_1AM_2^*\in \cl{B}.$ We proved that $\cl{M}\cl{A}\cl{M}^*$\\$\subset \cl{B}.$ Similarly we can prove that $\cl{M}^*\cl{B}\cl{M}\subset \cl{A}.  \qquad  \Box$

\section{TRO equivalence and spatial Morita equivalence}

The following definition is due to I. Todorov (personal communication).
\begin{definition}\label{4.1.d} Let $H_1,H_2$ be Hilbert spaces, $\cl{A}\subset B(H_1), \cl{B}\subset B(H_2)$ be $w^*$ closed algebras. If there exist linear spaces $\cl{U}\subset B(H_1,H_2), \cl{V}\subset B(H_2,H_1)$ such that $\cl{B}\cl{U}\cl{A}\subset \cl{U},$\;\; $\cl{A}\cl{V}\cl{B}\subset \cl{V},\;\; span( \cl{V}\cl{U})^{-w^*}=\cl{A}$ and $span( \cl{U}\cl{V})^{-w^*}=\cl{B}$  we say that the algebras $\cl{A},\cl{B}$ are \textbf{spatially Morita equivalent} and the system $(\cl{A},\cl{B},\cl{U},\cl{V})$ is a \textbf{spatial Morita context}.
\end{definition}

As we prove in Theorem \ref{4.4} and in remark \ref{4.5} spatial Morita equivalence is 
strictly weaker relation than TRO equivalence.

\begin{theorem}\label{4.1} Let $(\cl{A},\cl{B},\cl{U},\cl{V})$ be a spatial Morita context. Moreover we assume that the algebras $\cl{A},\cl{B}$ are unital. If $\phi =\mathrm {Map}(\cl{U})$ and $\psi =\mathrm {Map}(\cl{V})$ then

(i) $\cl{S}_{1,\phi }=\mathrm {Lat}(\cl{A}),\cl{S}_{2,\phi }=\mathrm {Lat}(\cl{B}),$ so the map $\phi :\mathrm {Lat}(\cl{A})\rightarrow \mathrm {Lat}(\cl{B})$ is a lattice isomorphism.

(ii) $\psi |_{\mathrm {Lat}(\cl{B})}=(\phi  |_{\mathrm {Lat}(\cl{A})})^{-1}$.
\end{theorem}
\textbf{Proof} Let $\zeta _1=\mathrm {Map}(\cl{A})$ and  $\zeta
_2=\mathrm {Map}(\cl{B}).$ Since $span( \cl{U}\cl{V})^{-w^*}=\cl{B}$ we get
$\zeta _2=\phi \circ \psi,$ hence $\zeta _2(pr(B(H_2)))\subset
\cl{S}_{2,\phi }$ or equivalently $\mathrm {Lat}(\cl{B})\subset
\cl{S}_{2,\phi }$.

Since $\zeta_1= \psi\circ \phi$, if
$P\in pr(B(H_1)),$ then
$$\cl{U}\cl{V}\phi (P)(H_2)\subset \cl{U}\psi
(\phi (P))(H_1)=\cl{U}\zeta _1(P)(H_1).$$
Also
$$\cl{U}\cl{A}P\subset \cl{U}P \Rightarrow
\cl{U}\zeta_1(P)(H_1)\subset \phi (P)(H_2) \Rightarrow
\cl{U}\cl{V}\phi (P)(H_2)\subset \phi (P)(H_2).$$
We proved that $\phi (P)^\bot \cl{U}\cl{V}\phi (P)=0$ and therefore
$\phi (P)^\bot \cl{B}\phi (P)=0$ for all $P\in pr(B(H_1)).$ It
follows that $\cl{S}_{2,\phi }\subset \mathrm {Lat}(\cl{B}),$ hence we have
equality.

Since $(\cl{A}^*,\cl{B}^*,\cl{U}^*,\cl{V}^*)$ is a spatial Morita context,
using the previous arguments we have $\cl{S}_{1,\phi }=\mathrm {Lat}(\cl{A}).$

Observe that $\phi :\mathrm {Lat}(\cl{A})\rightarrow \mathrm {Lat}(\cl{B})$ is a bijection which preserves order. Since $\psi \circ \phi =\zeta _1$, which is the identity on $\Lat{\cl A}$, it follows that $\psi\circ  (\phi|_{\mathrm {Lat}(\cl{A})})=Id|_{\mathrm {Lat}(\cl{A})}.$ Similarly $\phi \circ  (\psi |_{\mathrm {Lat}(\cl{B})})=Id|_{\mathrm {Lat}(\cl{B})}.$ The conclusion is that $\psi |_{\mathrm {Lat}(\cl{B})}= (\phi|_{\mathrm {Lat}(\cl{A})})^{-1}. \qquad \Box$

\begin{remark}\label{4.2}\em{ If $\cl{A},\cl{B}$ are spatially Morita equivalent unital algebras and the algebra $\cl{B}$ is reflexive, the algebra $\cl{A}$ is reflexive too. Indeed, let $(\cl{A},\cl{B},\cl{U},\cl{V})$ be a spatial Morita context. If $E,F$ are projections such that $E\cl{B}F=0,$ we have $E\cl{U}(\cl{VU})\cl{V}F=E\cl{(UV)(UV)}F=0$ hence $E\cl{UAV}F=0$ and so $E\cl{U}\Ref{\cl{A}}\cl{V}F=0.$ We proved that $\cl{U}\Ref{\cl{A}}\cl{V}\subset \cl{B}.$ But $$ \cl{U}\Ref{\cl{A}}\cl{V}\subset \cl{B} \Rightarrow \cl{V}\cl{U}\Ref{\cl{A}}\cl{V}\cl{U}\subset \cl{V}\cl{B}\cl{U}\subset \cl{A} \Rightarrow \cl{A}\Ref{\cl{A}}\cl{A}\subset \cl{A}.$$ Since the algebra $\cl{A}$ is unital we have $\Ref{\cl{A}}=\cl{A}.$

}
\end{remark}

\medskip

In the following theorem we prove that the converse of the above theorem is true for the case 
of CSL algebras.

\begin{theorem}\label{4.3} Two CSL algebras $\cl{A}$
and $\cl{B}$ are spatially Morita equivalent if and only if they have isomorphic lattices.
\end{theorem}
\textbf{Proof} By Theorem \ref{4.1} it suffices to show that a lattice isomorphism between CSL's
induces spatial Morita equivalence of the corresponding algebras. Suppose that $\cl{A}\subset B(H_1)$ and $\cl{B}\subset B(H_2).$ \\
Let $\cl{S}_1=\mathrm {Lat}(\cl{A}), \cl{S}_2=\mathrm {Lat}(\cl{B}),$ let
$\phi :\cl{S}_1\rightarrow \cl{S}_2$ be a lattice
isomorphism and
$$\cl{U}=\{T\in B(H_1,H_2): \phi (L)^\bot TL=0\:\text{for\;all\;}L\in \cl{S}_1\},$$
$$\cl{V}=\{S\in B(H_2,H_1): L^\bot S\phi (L)=0\:\text{for\;all\;}L\in \cl{S}_1\}.$$
 It is easily verified that $span( \cl{VU})$ is an ideal of $\cl A$.
Indeed if $V\in\cl V, U\in\cl U$ and $A\in \cl A$ then for all $L\in \cl S_1$ we have
$$AV\phi(L)=ALV\phi(L)=LALV\phi(L)$$ so $L^\bot AV\phi(L)=0$, hence $AV\in \cl V$.
Similarly $UA\in \cl U$, showing that $\cl{A}span( \cl{VU})\cl{A}\subset span( \cl{VU})$.
Also
\begin{align*}
VUL= V\phi(L)UL=LV\phi(L)UL
\end{align*}
and so $L^\bot VUL=0$; hence $VU\in \cl A$. It follows that $\Ref{\cl{VU}}\subset\cl{A}$. We shall prove that equality holds.

By Theorems 3.3, 4.4 in \cite{erd} we have that $$\cl{S}_{1,\mathrm {Map}(\cl{U})}=\cl{S}_1,\;\;\cl{S}_{2,\mathrm {Map}(\cl{U})}=\cl{S}_2,\;\; \mathrm {Map}(\cl{U})|_{\cl{S}_1}=\phi ,$$  $$\cl{S}_{1,\mathrm {Map}(\cl{V})}=\cl{S}_2,\;\; \cl{S}_{2,\mathrm {Map}(\cl{V})}=\cl{S}_1,\;\; \mathrm {Map}(\cl{V})|_{\cl{S}_2}=\phi^{-1}.$$
Let $\cl{W}=\Ref{\cl{V}\cl{U}}$ and $\zeta =\mathrm {Map}(\cl{W}).$ It follows that $\zeta =\mathrm {Map}(\cl{V})\circ \mathrm {Map}(\cl{U}).$ Also since $\cl{W}^*=\Ref{\cl{U}^*\cl{V}^*}$ we have $\zeta^* =\mathrm {Map}(\cl{U}^*)\circ \mathrm {Map}(\cl{V}^*),$ hence $\cl{S}_{2,\zeta ^*}\subset \cl{S}_{2,\mathrm {Map}(\cl{U}^*)}=(\cl{S}_{1,\mathrm {Map}(\cl{U})})^\bot =(\cl{S}_1)^\bot .$\\
We conclude that $\cl{S}_{1,\zeta }\subset \cl{S}_1.$
So if $L\in \cl{S}_{1,\zeta }$ we have  $\Map{\cl U}(L)\in \cl S_2$ and
$$\zeta (L)=\mathrm {Map}(\cl{V})\circ \mathrm {Map}(\cl{U})(L)=\phi^{-1}\circ \phi(L)=L,$$ 
hence

\begin{align*}\qquad \cl{W}=&\{T: \zeta (L)^\bot TL=0\;\text{for\;all\;}L\in \cl{S}_{1,\zeta }\}\\=
&\{T: L^\bot TL=0\;\text{for\;all\;}L\in \cl{S}_{1,\zeta }\}\\\supset 
&\{T: L^\bot TL=0\;\text{for\;all\;}L\in \cl{S}_{1}\}=\cl{A}
\end{align*}
Since $\cl{V}\cl{U}\subset \cl{A}$ we obtain the equality $\cl{A}=\Ref{\cl{V}\cl{U}}.$

Observe that the space $span( \cl{V}\cl{U})^{-w^*}$ is a masa bimodule, so it contains the space 
$\cl{A}_{min}.$ But the space $\cl{A}_{min}$ is a unital space and the space $span( \cl{V}\cl{U})^{-w^*}$ is an ideal of $\cl{A}.$ It follows that $\cl{A}=span( \cl{V}\cl{U})^{-w^*}.$ Similarly we can prove that $\cl{B}=span( \cl{U}\cl{V})^{-w^*}. \qquad \Box$
\begin{remark}
We do not know whether the `product' $span( \cl{V}\cl{U})^{-w^*}$ of two reflexive
spaces,or even reflexive masa bimodules, is necessarily reflexive.
\end{remark}
\begin{theorem}\label{4.4} If $\cl{A}, \cl{B}$ are $w^*$ closed TRO equivalent algebras 
then they are spatially Morita equivalent.
\end{theorem}
\textbf{Proof} Let $\cl{M}$ be an essential TRO such that $\cl{A}
\stackrel{\cl{M}}{\sim}\cl{B}$ and put
$\cl{A}_{\cl{M}}=span( \cl{A},\cl{M}^*\cl{M})^{-w^*},$ $
\cl{B}_{\cl{M}}=span( \cl{B},\cl{M}\cl{M}^*)^{-w^*}.$ We recall from
Proposition \ref{2.9} that $\cl{A}_{\cl{M}}
\stackrel{\cl{M}}{\sim}\cl{B}_{\cl{M}}$ and
$\cl{A},$\;(respectively $\cl{B}$) is an ideal of
$\cl{A}_{\cl{M}},$\;(respectively $\cl{B}_{\cl{M}}$).

 Let $\cl{U}=span( \cl{B}\cl{M})^{-w^*}$ and
$\cl{V}=span( \cl{M}^*\cl{B}_{\cl{M}})^{-w^*}.$ We shall show that the
system $(\cl{A},\cl{B},\cl{U},\cl{V})$ is a spatial Morita
context.

(i) Since $\cl{B}$ is an ideal of $\cl{B}_{\cl{M}}$ we have
$\cl{B}(\cl{BM})(\cl{M}^*\cl{BM})\subset \cl{BM}\subset \cl{U}.$
 Since $\cl{U}=span( \cl{B}\cl{M})^{-w^*}$ and
$\cl{A}=span( \cl{M}^*\cl{B}\cl{M})^{-w^*}$, this implies that
$\cl{B}\cl{U}\cl{A}\subset \cl{U}.$

(ii) Similarly, the relation  $(\cl{M}^*\cl{BM})(\cl{M}^*\cl{B}_{\cl{M}})\cl{B}\subset \cl{M}^*\cl{B}$ implies $\cl{A}\cl{V}\cl{B}\subset \cl{V}$ since  $\cl{V}=span( \cl{M}^*\cl{B}_{\cl{M}})^{-w^*}$.

(iii) Observe that $\cl{M}^*\subset \cl{V},$ hence
$\cl{M}^*\cl{B}\cl{M}\subset \cl{V}\cl{U}.$ It follows that\linebreak
$span( \cl{V}\cl{U})^{-w^*}\supset \cl{A}.$ Since $\cl{B}$ is an ideal of $\cl{B}_{\cl{M}}$ we have $\cl{M}^*\cl{B}_{\cl{M}}\cl{B}\cl{M}\subset \cl{M}^*\cl{B}\cl{M}\subset \cl{A}.$
But $\cl{V}=span( \cl{M}^*\cl{B}_{\cl{M}})^{-w^*}$ and $\cl{U}=span( \cl{B}\cl{M})^{-w^*}.$
It follows that $\cl{V}\cl{U}\subset \cl{A}.$ Therefore $\cl{A}=span( \cl{V}\cl{U})^{-w^*}.$

(iv) Since $\cl{B}$ is an ideal of $\cl{B}_{\cl{M}}$ we have
$\cl{B}\cl{M}\cl{M}^*\cl{B}_{\cl{M}}\subset \cl{B}$ hence
$\cl{U}\cl{V}\subset \cl{B}.$ Now, observe that
\begin{align*}
& \cl{M}\cl{M}^*\cl{B}_{\cl{M}}\subset \cl{M}\cl{V}\Rightarrow 
span( \cl{M}\cl{M}^*)^{-w^*}\cl{B}_{\cl{M}}\subset
span( \cl{M}\cl{V})^{-w^*}\Rightarrow \\ & \cl{B}_{\cl{M}}\subset
span( \cl{M}\cl{V})^{-w^*}.\end{align*} Also since $\cl{U}\cl{V}\supset
\cl{B}\cl{M}\cl{V},$ it follows that $$span( \cl{U}\cl{V})^{-w^*}\supset
\cl{B}span( \cl{M}\cl{V})^{w^*}\supset \cl{B}\cl{B}_{\cl{M}}=\cl{B}.$$
We proved that $\cl{B}=span( \cl{U}\cl{V})^{-w^*}.$ The proof is
complete.$\qquad  \Box$

\begin{remark}\label{4.5}
Spatial Morita equivalence does not imply TRO equivalence. There exist multiplicity free
nests $\cl S_1, \cl S_2$ which are isomorphic but the algebras $\cl S_1'', \cl S_2''$ are not isomorphic. For an example, see \cite[Example 7.19]{dav}.

Thus isomorphism of the lattices does not guarantee TRO equivalence, even for multiplicity free nest algebras.
\end{remark}

\begin{theorem}\label{4.6} Let $\cl{A}\subset B(H_1),\;\; \cl{B}\subset B(H_2)$
be separably acting CSL algebras with isomorphic lattices. Then $\cl A$ is
synthetic if and only if $\cl B$ is synthetic. In fact, if  $\phi: \Lat{\cl{A}}\to \Lat{\cl{B}}$ is a lattice isomorphism and $$\cl{U}=\{T\in B(H_1,H_2): \phi (L)^\bot TL=0\text{\;for\;all\;}L\in \Lat{\cl{A}}\}.$$ Then  $\cl{A}$ (and $\cl B$) is synthetic if and only if $\cl{U}$ is synthetic.
\end{theorem}
\textbf{Proof} Let $\cl S_1=\Lat{\cl A}$,
$$\cl{L}=\{\phi (L)\oplus L: L\in \cl{S}_1\}$$ 
and
$$\cl{V}=\{S\in B(H_2,H_1): L^\bot S\phi (L)=0\text{\;for\;all\;}L\in \cl{S}_1\}.$$
By Theorem \ref{4.3} we have that $$\cl{A}=span( \cl{V}\cl{U})^{-w^*} \text{\;\;and\;\;} \cl{B}=span( \cl{U}\cl{V})^{-w^*}.$$

It is shown in \cite[Proposition 4.2]{kt} that, if $\cl{C}=\Alg{\cl{L}}$,
$$\cl{C}=\left[\begin{array}{clr}\cl{B}&\cl{U}\\\cl{V} &\cl{A}\end{array}\right]\text{\;\;and\;\;}\cl{C}_{min}\subset
\left[\begin{array}{clr}\cl{B}_{min} & \cl{U}_{min} \\
\cl{V}_{min} & \cl{A}_{min}\end{array}\right].$$

We shall show that
\begin{equation}\label{41}
\cl{C}_{min}=
\left[\begin{array}{clr}\cl{B}_{min} & \cl{U}_{min} \\
\cl{V}_{min} & \cl{A}_{min}\end{array}\right].
\end{equation}

Indeed if $\cl{W}$ is any $w^*$ closed masa bimodule such
that $\Ref{\cl{W}}=\cl{C}$ and if $Q=0\oplus I$ then
$$\Ref{Q^\bot\cl{W}Q}=Q^\bot \cl{C}Q=
\left[\begin{array}{clr}0&\cl{U}\\0&0\end{array}\right].$$
It follows that
$$\left[\begin{array}{clr}0&\cl{U}_{min}\\0&0\end{array}\right]
=\left[\begin{array}{clr}0&\cl{U}\\0&0\end{array}\right]_{min}
        \subset Q^\bot \cl{W}Q\subset \cl{W}.$$
Now taking $\cl{W}=\cl{C}_{min}$ we obtain
$$\left[\begin{array}{clr}0&\cl{U}_{min}\\0&0\end{array}\right]\subset
\cl{C}_{min}.$$
 Similarly we can prove that $$\left[\begin{array}{clr}0& 0 \\ 0 & \cl{A}_{min} \end{array}\right], \left[\begin{array}{clr}0&0\\\cl{V}_{min} & 0\end{array}\right], \left[\begin{array}{clr} \cl{B}_{min} & 0 \\0&0\end{array}\right]\subset
\cl{C}_{min}$$
and (\ref{41}) follows.

If $E,F$ are projections such that $E\cl{U}_{min}\cl{V}F=0$ then
$E\cl{U}\cl{V}F=0$ hence $E\cl{B}F=0.$ It follows that
$\cl{B}\subset \Ref{\cl{U}_{min}\cl{V}}$ and therefore
$\cl{B}= \Ref{\cl{U}_{min}\cl{V}}.$ Similarly
$\cl{A}= \Ref{\cl{V}_{min}\cl{U}}.$

Now suppose that the algebra $\cl{A}$ is synthetic. Since $\cl{B}= \Ref{\cl{U}_{min}\cl{V}}$ we have that $\cl{B}_{min}\subset span( \cl{U}_{min}\cl{V})^{-w^*}$ so \begin{equation}\label{U}\cl{U}\subset \cl{B}_{min}\cl{U}\subset span( \cl{U}_{min}\cl{V}\cl{U})^{-w^*}\subset span( \cl{U}_{min}\cl{A})^{-w^*}=span( \cl{U}_{min}\cl{A}_{min})^{-w^*}.\end{equation}

Using (\ref{41}) we have  $$\left[\begin{array}{clr}0&\cl{U}_{min}\cl{A}_{min}\\0&0\end{array}\right]=\left[\begin{array}{clr}0&\cl{U}_{min}\\0&0\end{array}\right]\left[\begin{array}{clr}0&0\\0&\cl{A}_{min}\end{array}\right]\subset
\cl{C}_{min}.$$

It follows that $\cl{U}_{min}\cl{A}_{min}\subset \cl{U}_{min},$ hence $\cl{U}\subset \cl{U}_{min}$ (from equation (\ref{U})) and so the bimodule $\cl{U}$ is synthetic.

For the opposite direction we suppose that the bimodule $\cl{U}$ is synthetic. Since
$$\left[\begin{array}{clr}0&0\\0&\cl{V}_{min}\cl{U}_{min}\end{array}\right]=
\left[\begin{array}{clr}0&0\\\cl{V}_{min}&0\end{array}\right]
\left[\begin{array}{clr}0&\cl{U}_{min}\\0&0\end{array}\right],$$
again using (\ref{41}) we conclude that $\cl{V}_{min}\cl{U}_{min}\subset \cl{A}_{min}$
and therefore $\cl{V}_{min}\cl{U}\subset \cl{A}_{min}.$
Since $\cl{U}\cl{A}\subset \cl{U}$ it follows
that $\cl{V}_{min}\cl{U}\cl{A}\subset \cl{A}_{min}.$
But also  $\cl{A}_{min}\subset span( \cl{V}_{min}\cl{U})^{-w^*}$ since
$\cl{A}=\Ref{\cl{V}_{min}\cl{U}}$, and hence
$\cl{A}_{min}\cl{A}\subset \cl{A}_{min}$; therefore $\cl{A}\subset \cl{A}_{min}$ since
$\cl A_{min}$ is unital.

We have proved that the algebra $\cl{A}$ is synthetic if and only if the bimodule $\cl{U}$ is synthetic. Similarly one shows that  $\cl{U}$ is synthetic if and only if the algebra $\cl{B}$ is synthetic. $\qquad  \Box$

\section{TRO equivalence and CSL  algebras}

In this section we assume that all Hilbert spaces are separable. 
Thus the $w^*$ topology on bounded sets of operators is metrisable. We are going 
present some results on TRO equivalence of CSL algberas.

\begin{definition}\label{5.1.d} If $\cl{S}$ is a CSL and $L\in \cl{S}$ we denote by $L_\flat$
the projection $\vee \{M\in \cl{S}: M<L\}.$ Whenever $L_\flat <L$ we call
the projection $L-L_\flat $ an \textbf{atom} of $\cl{S}.$ If the CSL $\cl{S}$ has no atoms 
we say that it is a \textbf{continuous} CSL. If the identity operator is the sum of the atoms we say that 
 $\cl S$ is \textbf{totally atomic}.
\end{definition}

\begin{lemma}\label{5.2} Let $H_i, i=1,2$ be Hilbert spaces, $\cl{S}_i, i=1,2$ be commutative 
lattices (not necessarily complete) containing zero and the identity and let 
$\theta :\cl{S}_1\rightarrow \cl{S}_2$ be a lattice isomorphism. Then the map 
$\theta $ extends to a $*$ isomorphism $\rho :span( \cl{S}_1)^{-\|\cdot\|}\rightarrow span( \cl{S}_2)^{-\|\cdot\|}.$
\end{lemma}
\textbf{Proof} Using induction we shall prove that if $P_1,...P_n$ are projections of  $\cl{S}_1$
such that $\sum_{i=1}^nc_iP_i=0$ where $c_i\neq0$ for  $1\leq i\leq n$ then  $\sum_{i=1}^nc_i\theta (P_i)=0.$

 The claim clearly holds for $n=1,2$. Assume that it holds for $k\in \{1,...,n-1\}.$

Let $\sum_{i=1}^nc_iP_i=0$ where $c_i\neq0$ for  $1\leq i\leq n$ and put $A=\sum_{i=1}^nc_i\theta (P_i).$ It suffices to show that $\theta (P_k)A=0$ for all $k\in \{1,...n\}.$

Let $B=\theta (P_n)A.$ We shall show that $B=0.$

Multiply the equation  $\sum_{i=1}^nc_iP_i=0$ with $P_1\wedge P_n.$ This gives
$$(c_1+c_n)(P_1\wedge P_n)+c_2(P_2\wedge P_1\wedge P_n)+\ldots
+c_{n-1}(P_{n-1}\wedge P_1\wedge P_n)=0.$$
By the inductive hypothesis we have
$$(c_1+c_n)\theta (P_1\wedge P_n)+c_2\theta (P_2\wedge P_1\wedge P_n)+\ldots
+c_{n-1}\theta (P_{n-1}\wedge P_1\wedge P_n)=0,$$
hence $\theta (P_1)B=0.$

Similarly we can prove that $\theta (P_i)B=0$ for $1\leq i\leq n-1.$\\
Since $P_n=(-c_n)^{-1}\sum_{i\neq n}c_i P_i$ it follows that $P_n\leq \vee _{i\neq n}P_i$ hence $\theta (P_n)\leq \vee _{i\neq n}\theta (P_i),$ so $\theta (P_n)B=0$ and therefore $B=0.$ We proved that $\theta (P_n)A=0.$\\
Using the same method we have  $\theta (P_k)A=0, k=1,...n.$ So the claim holds.

The conclusion is that the map $$\rho :span( \cl{S}_1)\rightarrow
span( \cl{S}_2): \rho \left(\sum_{i=1}^nc_iP_i\right)=\sum_{i=1}^nc_i\theta
(P_i)$$ is well defined, and it is clearly a $*$ isomorphism.

We shall show that $\rho $ is norm continuous. {Let
$T=\sum_{i=1}^nc_iP_i\in span( \cl{S}_1)$ and let $c$ be in the
spectrum $\sigma (\rho (T))$ of $\rho (T)$.} Let
$\cl{S}_0$ be the smallest lattice containing the set
$\{0,P_1,...P_n,I\}.$ Then the space $span( \cl{S}_0)$ is a 
$C^*$-algebra which is contained in $span( \cl{S}_1).$
If $c$ is not in the spectrum $\sigma (T)$ of $T$,
the operator $S=(cI-T)^{-1}$ is contained in $span( \cl{S}_0)$ and hence in
$span( \cl{S}_1).$ Since $S(cI-T)=(cI-T)S=I$ we have $\rho (S)(cI-\rho
(T))=(cI-\rho (T))\rho (S)=I$, contradicting $c\in\sigma (\rho (T))$.

We proved that $\sigma (\rho (T))\subset\sigma (T).$ Therefore
$\|\rho (T)\|\leq \|T\|.$

We conclude that the map $\rho $
extends to a $*$ isomorphism from the $C^*$-algebra
$span( \cl{S}_1)^{-\|\cdot\|}$ onto $span( \cl{S}_2)^{-\|\cdot\|}.$
$\qquad \Box$

\begin{lemma}\label{5.3} Let $\cl{S}_1, \cl{S}_2$ be CSL's,
$\phi : \cl{S}_1\rightarrow \cl{S}_2$ be a lattice isomorphism, $P$ be the sum of the
atoms of $\cl{S}_1$ and $Q$ be the sum of the
atoms of $\cl{S}_2.$ Then there exists a $*$ isomorphism
$$\rho :\cl{S}_1''|_P\rightarrow \cl{S}_2''|_Q$$
such that $\rho (L|_P)=\phi (L)|_Q$ for all $L\in \cl{S}_1.$
\end{lemma}
\textbf{Proof} For any atom $A$ in a CSL $\cl S,$ we have 
$$A=\wedge \{L\in \cl S: AL=A\}-\vee \{L\in \cl S: AL=0\}.$$ Letting 
$\cl A_i$ be the set of atoms of $\cl S_i,$ we then find that the lattice isomorphism 
$\phi : \cl S_1\rightarrow \cl S_2$ induces a bijection $\rho _\alpha : \cl A_1\rightarrow 
\cl A_2,$ given by 
 $$\rho _\alpha (A)=\wedge \{\phi (L): L\in \cl S: AL=A\}-\vee \{\phi (L): L\in \cl S: AL=0\}.$$ 
Now $\cl{S}_1''|_P$ and $\cl{S}_2''|_Q$ are $*$ isomorphic to $l^\infty (\cl A_1)$ and 
$l^\infty (\cl A_2)$ respectively; identify  $\cl{S}_1''|_P$ with $l^\infty (\cl A_1)$ 
and likewise identify $\cl{S}_2''|_Q$ with $l^\infty (\cl A_2).$ The 
map $\rho_ \alpha $ extends to an isomorphism $\rho : l^\infty (\cl A_1)\rightarrow 
l^\infty (\cl A_2),$ which satisfies our requirements, because any projection in $l^\infty (\cl A_i)$ 
is the sum of the atoms it contains. $\qquad \Box$

\medskip

The following theorem is consequence of the above lemma.

\begin{theorem}\label{theorema1} Let $\cl S_1, \cl S_2$ be totally atomic CSL's. The algebras 
$\mathrm{Alg}(\cl S_1), \mathrm{Alg}(\cl S_2)$ are TRO equivalent if and only if the CSL's $\cl S_1,
 \cl S_2$ are isomorphic. 
\end{theorem}
\textbf{Proof} If $ \mathrm{Alg}(\cl S_1) $ and $ \mathrm{Alg}(\cl S_2)$ are TRO equivalent, 
then $\cl S_1$ and $\cl S_2$ are isomorphic by Theorem \ref{3.3}. Conversely, by the above 
lemma any lattice isomorphism $\phi : \cl S_1\rightarrow \cl S_2$ extends to a $*$ isomorphism 
from $\cl S_1^{\prime \prime}=(\Delta (\mathrm{Alg}(\cl S_1)
))^\prime $ onto $\cl S_2^{\prime \prime}=(\Delta (\mathrm{Alg}(\cl S_2)))^\prime .$ 
Using again Theorem \ref{3.3} we conclude that the algebras $ \mathrm{Alg}(\cl S_1) $ 
and $ \mathrm{Alg}(\cl S_2)$ are TRO equivalent. $\qquad \Box$

\medskip

For the general case of CSL algebras we present the following theorem.

\begin{theorem}\label{5.7.b} Let  $\cl{S}_1, \cl{S}_2$ be CSL's acting on Hilbert spaces
$H_1, H_2$ respectively, $P$ the sum of the atoms of $\cl{S}_1,$ $Q$ the sum of the atoms of
$\cl{S}_2$ and $\cl{A}=\mathrm {Alg}(\cl{S}_1),$ $\cl{B}=\mathrm {Alg}(\cl{S}_2),$ 
$\cl{A}_0=span( \cl{S}
_1)^{-\|\cdot\|},$ $\cl{B}_0=span( \cl{S}_2)^{-\|\cdot\|}.$ The following are equivalent:

 (i)  The algebras $\cl{A},\cl{B}$ are TRO equivalent.

(ii)There exists a lattice isomorphism $\phi : \cl{S}_1\rightarrow \cl{S}_2$
whose extension (Lemma \ref{5.2}) $\overline {\phi}  :\cl{A}_0 \rightarrow \cl{B}_0$ is
$w^*$-bicontinuous on the unit balls of $\cl{A}_0, \cl{B}_0.$

(iii) There exists a lattice isomorphism $\phi : \cl{S}_1\rightarrow \cl{S}_2$ such that  
if $\cl{L}=
\{L\oplus \phi (L): L\in \cl{S}_1\}$ then $$ \cl{L}^{\prime \prime}\cap (0\oplus \cl{B}_0^
{\prime \prime})=0,\;\;  \cl{L}^{\prime \prime}\cap (\cl{A}_0^{\prime \prime}\oplus 0)=0.$$

(iv) There exists a lattice isomorphism $\phi : \cl{S}_1\rightarrow \cl{S}_2$ such that  
if $\cl{L}$
 is as in (iii) then $$ \cl{L}^{\prime \prime}\cap (0\oplus \cl{B}_0^
{\prime \prime}Q^\perp )=0,\;\;  \cl{L}^{\prime \prime}\cap (\cl{A}_0^{\prime \prime}P^\perp 
\oplus 0)=0.$$

 Moreover if these conditions hold and
$$\Delta (\phi )=\{T\in B(H_1, H_2): TL=\phi (L)T \;\text{for\;all\;}L\in \cl{S}_1\}$$
then  $\cl{A} \stackrel{\Delta (\phi )}{\sim}\cl{B}.$
\end{theorem}
\textbf{Proof (i)$\Rightarrow $(ii)}\\
This is obvious by Theorem \ref{3.3}
 \medskip

\noindent\textbf{(ii)$\Rightarrow $(i)}\\
Suppose that $\phi : \cl{S}_1\rightarrow \cl{S}_2$ is a lattice isomorphism
whose extension by Lemma \ref{5.2},
$\overline {\phi}  :\cl{A}_0 \rightarrow \cl{B}_0$ is $w^*$-bicontinuous on the unit balls.
 By \cite[Lemma 10.1.10]{kr} the map $\overline{\phi }$
(respectively \;$\overline{\phi }^{-1}$ ) extends to a $w^*$-continuous
homomorphism from $(\cl{S}_1)''$ to $(\cl{S}_2)''$ (respectively from
$(\cl{S}_2)''$ to $(\cl{S}_1)''$). One can check that the extensions are mutual
inverses. (The assumption that the map $\overline{\phi }$ is $w^*$-continuous
doesn't guarantee that its inverse is $w^*$-continuous. See exercise 10.5.30 in \cite{kr}).

 Now Theorem \ref{3.3} shows that the algebras $\cl{A}$ and $\cl{B}$ are TRO equivalent.

\medskip

\noindent\textbf{(i)$\Rightarrow $(iii)}\\
If the algebras $\cl{A}, \cl{B}$ are TRO equivalent, by Theorem \ref{3.3} there exists a 
lattice isomorphism
$\phi : \cl{S}_1 \rightarrow  \cl{S}_2$ which extends to a $*$ isomorphism
$\rho : \cl{S}_1^{\prime \prime} \rightarrow  \cl{S}_2^{\prime \prime}.$ We can verify that
$\cl{L}^{\prime \prime}=\{A\oplus \rho (A): A\in \cl{A}_0^{\prime \prime}\}$ hence  
$$ \cl{L}^{\prime \prime}\cap (0\oplus \cl{B}_0^
{\prime \prime})=0,  \cl{L}^{\prime \prime}\cap (\cl{A}_0^{\prime \prime}\oplus 0)=0.$$

\noindent\textbf{(iii)$\Rightarrow $(iv)}\\
This is obvious.

\medskip

\noindent\textbf{(iv)$\Rightarrow $(i)}\\
It suffices to show that $\phi $ extends to a $*$ isomorphism from $\cl{S}_1^{\prime \prime}$
onto $\cl{S}_2^{\prime \prime}.$ If it does not by (ii) one of the maps
 $\overline {\phi}  :\cl{A}_0 \rightarrow \cl{B}_0,$
$\overline {\phi}^{-1}  :\cl{B}_0 \rightarrow \cl{A}_0$ will not be $w^*$ continuous 
on the unit ball.
Suppose that $ \overline{\phi } $ is not $w^*$ continuous on the unit ball. Then there
exists a net $(A_i)\subset Ball(\cl{A}_0)$ which converges in the $w^*$ topology to $0$
while the net $(\overline{\phi }(A_i))$ converges to a nonzero operator $B\in \cl{B}_0^{\prime
\prime}.$ Since the restriction of $\phi $ on the lattice $\cl S_1|_P$ extends (Lemma \ref{5.3}) 
to a $*$ isomorphism from
$\cl{A}_0^{\prime \prime}|_P$ onto $\cl{B}_0^{\prime \prime}|_Q$ and the net $(A_iP)$
 converges to $0$ the net $(\overline{\phi }(A_i)Q)$ converges to $0$ too. Therefore $BQ=0.$

 Observe that $(L-L_\flat)\oplus (\phi(L)-\phi(L)_\flat)\in \cl{L}^{\prime \prime}$ for all
$L\in \cl{S}_1,$ hence $P\oplus Q \in \cl{L}^{\prime \prime}.$ It follows that $A_iP^\perp \oplus
\overline{\phi }(A_i)Q^\perp  \in \cl{L}^{\prime \prime}$ for every index $i$ and so $0\oplus
BQ^\perp \in \cl{L}^{\prime \prime}.$ This is a contradiction because $BQ^\perp \neq 0.$
 The proof for the case where $ \overline{\phi }^{-1} $ is not $w^*$ continuous on the unit
ball is similar.

\medskip

 Now suppose that conditions (ii) to (v) hold and let $\rho : \cl{S}_1^{\prime
\prime} \rightarrow \cl{S}_2^{\prime \prime}$ be the extension of $\phi .$
 By Lemma \ref{3.1} the
 space $$\cl M=\{T\in B(H_1,H_2): TA=\rho (A)T\;
\text{for\;all\;}A\in \cl{S}_2^{\prime\prime}\}$$
is an essential TRO. Since the space $\Delta (\phi )$
contains $\cl M$ it is essential too.
We can easily verify that
$\Delta (\phi )^*\cl{B}\Delta (\phi) \subset \cl{A}$ and
$\Delta (\phi)\cl{A}\Delta (\phi)^*\subset \cl{B}.$
By Proposition \ref{2.1} we have
$\cl{A} \stackrel{\Delta (\phi )}{\sim}\cl{B}.$ $\qquad \Box$

\begin{remark}\label{5.8}
By Theorem \ref{3.3}, if the algebras $\cl{A}$ and $\cl{B}$ are TRO equivalent and
$\cl{A}$ is a CSL algebra then so is $\cl{B}$.
\end{remark}

\medskip

In the special case of nest algebras we have the following result:.

\begin{theorem}\label{theorema2} All nest algebras with continuous nests are TRO equivalent. 
\end{theorem}
\textbf{Proof} If $\cl X$ is a subset of some $B(H)$ we denote by $\cl X^\infty $ 
the set of all operators of the form $X^\infty =X\oplus X \oplus....,$ 
where $X\in \cl X,$ acting on $B(H^\infty  ).$ If $\cl S_1, \cl S_2$ are continuous 
nests, the nests $\cl S_1^\infty , \cl S_2^\infty $ are also continuous and of 
multiplicity $\infty .$ It follows from \cite[Theorem 7.24]{dav} that the nests 
$\cl S_1^\infty , \cl S_2^\infty $ are unitarily equivalent. So there exists a $*$ isomorphism 
from $ (\cl S_1^{\prime \prime})^\infty $ onto $(\cl S_2^{\prime \prime})^\infty $ 
 mapping $\cl S_1^\infty$ onto  $\cl S_2^\infty.$ Now taking compositions with the maps 
$ \cl S_i^{\prime \prime} \rightarrow  (\cl S_i^{\prime \prime})^\infty, X\rightarrow X^\infty ,
 i=1,2 $ we obtain a  $*$ isomorphism from  $\cl S_1^{\prime \prime}$ onto 
$  \cl S_2^{\prime \prime}$ mapping $\cl S_1$ onto $\cl S_2.$ The conclusion comes 
from Theorem \ref{3.3}. $\qquad \Box$

\bigskip

{\em Acknowledgement:} I would like to express appreciation to
Prof.~A.~Katavolos for his helpful comments and suggestions during
the preparation of this work, which is part of my doctoral thesis.
I am also indebted to
Dr.~I.~Todorov for helpful discussions and important suggestions. Finally, I am grateful  
to the referee for carefully reading the paper and making significant improvements. In 
particular I wish to thank him for pointing out an error in a previous version and for 
providing a quicker 
proof of Lemma \ref{5.3}.

\end{document}